\documentclass{article}
\usepackage{amsfonts}
\usepackage{amsmath,amsthm}

\usepackage[all]{xy}
\usepackage{graphicx}

\usepackage{amssymb}
\usepackage{latexsym}

\DeclareMathAlphabet\EuFrak{U}{euf}{m}{n}	
\SetMathAlphabet\EuFrak{bold}{U}{euf}{b}{n}	

\newcommand{\ra}{\rightarrow}

\newcommand{\hra}{\hookrightarrow}

\newcommand{\ovl}{\overline}

\newcommand{\wa}{\widehat}

\newcommand{\sC}{{\it C*}-}
\newcommand{\bC}{{\mathbb C}}

\newcommand{\bE}{{\mathbb E}}

\newcommand{\bT}{{\mathbb T}}

\newcommand{\bZ}{{\mathbb Z}}

\newcommand{\bN}{{\mathbb N}}

\newcommand{\mA}{\mathcal A}
\newcommand{\mB}{\mathcal B}
\newcommand{\mC}{\mathcal C}

\newcommand{\mE}{\mathcal E}
\newcommand{\mF}{\mathcal F}

\newcommand{\mK}{\mathcal K}
\newcommand{\mL}{\mathcal L}
\newcommand{\mM}{\mathcal M}
\newcommand{\mN}{\mathcal N}
\newcommand{\mO}{\mathcal O}

\newcommand{\mR}{\mathcal R}
\newcommand{\mT}{\mathcal T}

\newcommand{\mU}{\mathcal U}

\newcommand{\wE}{\wa {\mathcal E}}

\newcommand{\pic}{ { {\mathrm {Pic}} (X ; \mA)} }

\newcommand{\coe}{\mO_{\mE}}

\newcommand{\com}{\mO_{\mM}}

\newtheorem{thm}{Theorem}[section]
\newtheorem{cor}[thm]{Corollary}
\newtheorem{lem}[thm]{Lemma}
\newtheorem{prop}[thm]{Proposition}

\newtheorem{defn}[thm]{Definition}

\theoremstyle{definition}
\newtheorem{ex}{Example}[section]

\theoremstyle{remark}
\newtheorem{rem}{Remark}[section]

\numberwithin{equation}{section}

\begin{document}

\author{{\sf Ezio Vasselli}
                         \\{\it Dipartimento di Matematica}
                         \\{\it Universit\`a La Sapienza di Roma}
			 \\{\it P.le Aldo Moro, 2 - 00185 Roma - Italy }
 			 \\{\it (c/o Sergio Doplicher) }
                         \\{\sf vasselli@mat.uniroma2.it}}

\title{ Bundles of \sC algebras \\and\\the $KK(X ; - , - )$-bifunctor}
\maketitle

\begin{abstract}

An overview about {\it C*}-algebra bundles with a $\mathbb{Z}$-grading is presented, with particular emphasis on classification questions. In particular, we discuss the role of the representable $KK(X ; - , - )$-bifunctor introduced by Kasparov. As an application, we consider Cuntz-Pimsner algebras associated with vector bundles, and give a classification in terms of $K$-theoretical invariants in the case in which the base space is an $n$-sphere.

\bigskip

\noindent {\bf AMS Subj. Class.:} 46L05, 46L80.


\end{abstract}

\section{Introduction.}
\label{intro}

The classification of \sC algebras by $K$-theoretical invariants is a rich and interesting topic; the relative program particularly succeeded in the case of simple, nuclear, purely infinite \sC algebras (\cite{Kir03,Phi03}).

In order to extend such results to the case of non-simple \sC algebras, it is natural to consider a particular class, namely \sC algebra bundles over a locally compact Hausdorff space $X$. In order to find good invariants, in this case the better-behaved tool is the representable $KK ( X ; - , - )$-theory introduced by Kasparov in \cite{Kas88}, which takes into account the bundle structure of the underlying \sC algebra. $KK ( X ; - , - )$-theory has been recently extended to the case in which $X$ is a $T_0$ space, in order to consider primitive ideal spectra of \sC algebras (\cite{Kir00}).

Aim of the present paper is to present an overview about \sC algebra bundles with a $\bZ$-grading, and their associated $KK ( X ; - , - )$-theoretical invariants. Our main motivation arises from the case of the universal \sC algebra of a vector bundle $\mE \ra X$, which is constructed as the Cuntz-Pimsner algebra associated with the bimodule of continuous sections of $\mE$. Such a \sC algebra has a natural structure of a $\bZ$-graded bundle over $X$, with fibre the well-known Cuntz algebra. We are interested to classify such \sC algebras in terms of properties of the underlying vector bundles.

All the material presented in the present paper appeared elsewhere (in some different form), with the exception of our main result Thm.\ref{prop_sn}, where we classify \sC algebras of vector bundles in the case in which the base space is an $n$-sphere.

It is aim of the present work to be self-contained: the reader is assumed to be familiar at an elementary level with \sC algebra theory (\cite{Ped}), and $K$-theory (\cite{Ati,Bla}). In the case of results proved elsewhere, the proofs will be sketched or omitted. Some of the results exposed in the present paper appear in \cite{Vas04}.

\section{Bundles and $C_0(X)$-algebras.}
\label{fields}

Let $X$ be a locally compact Hausdorff space, $C_0(X)$ the \sC algebra of complex-valued, continuous, vanishing at infinity functions on $X$. A {\em continuous bundle} of \sC algebras over $X$ is a \sC algebra $\mF$, equipped with a faithful family of epimorphisms $\left\{ \pi_x : \mF \ra \mF_x \right\}_{x \in X}$ such that, for every $a \in \mF$, the {\em norm function} $\left\{ x \mapsto \left\| \pi_x (a) \right\|  \right\}$ belongs to $C_0(X)$; furthermore, $\mF$ is required to be a nondegenerate $C_0(X)$-module w.r.t. pointwise multiplication $f , a \mapsto \left\{ f(x) \cdot \pi_x (a) \right\}$, $f \in C_0(X)$. If $X$ is compact, we consider the analogous notion by using the \sC algebra $C(X)$ of continous functions on $X$.

\begin{ex}
Let $\mA$ be a \sC algebra. Then, the \sC algebra tensor product $C_0(X) \otimes \mA$ is a continuous bundle, called the {\bf trivial bundle}. To be more concise, we define
\begin{equation}
\label{def_xa}
X \mA := C_0(X) \otimes \mA \ .
\end{equation} 
\end{ex}

The above notion of continuous bundle has been given in \cite{KW95}: it is a simplified version of the classical notion of {\em continuous field} (see \cite[\S 10]{Dix}). We refer to the last-cited reference for the notions of {\em restriction} (\cite[10.1.7]{Dix} and {\em local triviality} (\cite[10.1.8]{Dix}), which are the analogues to well-known notions in the setting of topological bundles.

Let $\mA$ be a \sC algebra. To be more concise, we will call $\mA${\em - bundle} a locally trivial continuous bundle $\mF$ with fibre $\mF_x \equiv \mA$, $x \in X$.

A $C_0(X)$-{\em algebra} is a \sC algebra $\mA$, equipped with a nondegenerate morphism from $C_0(X)$ into the centre of the multiplier algebra $M(\mA)$ (\cite[\S 2]{Kas88}); in the sequel, we will identify $C_0(X)$ with the image in $M(\mA)$. \sC algebra morphisms commuting with the $C_0(X)$-actions are called {\em $C_0(X)$-morphisms}. We denote by ${\bf aut}_X \mA$ the group of $C_0(X)$-automorphisms of $\mA$. It is proved in \cite{Nil96} that the category of $C_0(X)$-algebras is equivalent to the one of 'upper semicontinuous bundles'; the fibre of $\mA$ over $x \in X$ is defined as follows: we consider the closed ideal $I_x := C_0(X - \left\{ x \right\}) \cdot \mA \subset \mA$, and define $\mA_x := \left. \mA \right/ I_x$. In particular, every continuous bundle is a $C_0(X)$-algebra. We will denote by $\otimes_X$ the minimal tensor product with coefficients in $C_0(X)$ (\cite{Bla95}).

\section{Hilbert bimodules and Cuntz-Pimsner algebras.}

For basic notions and terminology about {\em Hilbert bimodules}, we refer to \cite[\S 13]{Bla}.

Let $\mB$ be a \sC algebra, $\mM$ a right Hilbert $\mB$-module. We denote by $L(\mM)$ the \sC algebra of (bounded) adjointable right $\mB$-module operators on $\mM$, and by $K (\mM)$ the ideal of compact right $\mB$-module operators of the type
\begin{equation}
\label{def_theta}
\theta_{\psi , \psi'} \varphi := 
\psi \cdot \left \langle \psi' , \varphi \right \rangle \ \ ,
\end{equation}
\noindent where $\psi , \psi' , \varphi \in \mM$ and $\left \langle \cdot , \cdot \right \rangle$ denotes the $\mB$-valued scalar product.

Let $\mM$ be a Hilbert $\mA$-$\mB$-bimodule. In the sequel, if not specified, we will identify elements of $\mA$ with the corresponding operators in $L(\mM)$, by assuming that the morphism $\mA \ra L(\mM)$ is injective.

\begin{defn}
\label{2morphism}
Let $\mA$, $\mB$ be a \sC algebras, $\mM$ a Hilbert $\mA$-bimodule, $\mN$ a Hilbert $\mB$-bimodule. A {\bf covariant morphism} from $\mM$ into $\mN$ is a pair $(\beta , \eta)$, where $\beta : \mM \ra \mN$ is a Banach space map, $\eta : \mA \ra \mB$ is a \sC algebra morphism, and the following properties are satisfied for $a \in \mA$, $\psi , \psi' \in \mM$:
\[  
\beta (a \psi) = \eta (a) \beta (\psi) \ \ , \ \ 
\beta (\psi a) = \beta (\psi) \eta (a) \ \ , \ \
\left \langle \beta (\psi) , \beta(\psi') \right \rangle  = 
      \eta \left \langle \psi , \psi' \right \rangle \ ,
\]
\noindent where $\left \langle \cdot , \cdot \right \rangle$ denotes the $\mA$-valued (resp. $\mB$-valued) scalar product of $\mM$ (resp. $\mN$).
\end{defn}

\begin{ex}
\label{ex_pullback}
Let $\alpha : \mA \ra \mB$ be a \sC algebra isomorphism, $\mM$ a Hilbert $\mA$-bimodule. We introduce a Hilbert $\mB$-bimodule $\mM_\alpha$, defined as the set $\mM_\alpha := \left\{ \underline \psi , \psi \in \mM \right\}$ endowed with the vector space structure induced by $\mM$. The Hilbert $\mB$-bimodule structure is defined as follows:
\[
b \underline \psi := \underline {\alpha^{-1} (b) \psi} \ \ , \ \
\underline \psi b := \underline {\psi \alpha^{-1} (b)} \ \ , \ \
\left \langle \underline \psi , \underline \psi' \right \rangle := 
       \alpha \left \langle \psi , \psi' \right \rangle \ .
\]
\noindent We call $\mM_\alpha$ the {\em pullback bimodule} of $\mM$. Let now $\beta (\psi) := \underline \psi$, $\psi \in \mM$; it is clear that the pair $(\beta , \alpha)$ is a covariant isomorphism from $\mM$ onto $\mM_\alpha$. Viceversa, if $(\beta , \alpha)$ is a covariant isomorphism from a Hilbert $\mA$-bimodule $\mM$ into a Hilbert $\mB$-bimodule $\mN$, then $\mM_\alpha$ is isomorphic to $\mN$ as a Hilbert $\mB$-bimodule. 
\end{ex}

%
%

Let $\mA$ be a \sC algebra, $\mM$ a Hilbert $\mA$-bimodule. The {\em Cuntz-Pimsner C*-algebra} ({\em CP-algebra}, in the sequel) associated with $\mM$ has been introduced in \cite{Pim93}; it is obtained by a universal construction, and supplies a generalization of crossed products by $\bZ$ (see Ex.\ref{ex_mor} below) and the well-known Cuntz algebras $\mO_d$, $d \in \bN$ (\cite{Cun77}). We will denote by $\com$ the CP-algebra associated with $\mM$.

In order to simplify the exposition, we give a description of $\com$ in the case in which $\mA$ has identity $1$ and $\mM$ is finitely generated as a right Hilbert $\mA$-module. Let $\left\{ \psi_l \right\}_{l=1}^n \subset \mM$ be a finite set of generators, $\left \langle \cdot , \cdot \right \rangle$ the $\mA$-valued scalar product; then, for every index $l$, $a \in \mA$, we find $a \psi_l = \sum_m \psi_m a_{ml}$, $a_{ml} := \left \langle \psi_m , a \psi_l \right \rangle \in \mA$. We consider the universal $*$-algebra $^0 \com$ generated by $\left\{ \psi_l \right\}$, $\mA$, with relations
\begin{equation}
\label{eq_rcom}
\psi_l^* \psi_m = \left \langle \psi_l , \psi_m \right \rangle 
\ \ , \ \
a \psi_l = \sum_m \psi_m a_{ml}
\ \ , \ \
\sum_l \psi_l \psi_l^* = 1 \ .
\end{equation}
\noindent Note that every $\psi \in \mM$ appears as an element of $^0 \com$, in fact $\psi = \sum_l \psi_l (\psi_l^* \psi)$. It can be proved that there exists a {\em unique} (faithful) \sC norm on $^0 \com$ such that the {\em circle action}
\begin{equation}
\label{def_ca}
\alpha_z (\psi) := z \psi \ \ , \ \ z \in \bT , \psi \in \mM
\end{equation}
\noindent extends to an (isometric) automorphic action. The resulting \sC algebra is the CP-algebra $\com$, naturally endowed with the action $\alpha : \bT \ra {\bf aut}\com$. We introduce the notation
\begin{equation}
\label{def_comk}
\com^k := \left\{ t \in \com : \alpha_z (t) = z^k t \right\} 
\ \ , \ \
k \in \bZ \ .
\end{equation}

\begin{ex}
The Cuntz algebra $\mO_d$, $d \in \bN$, is obtained in the case $\mA = \bC$, $\mM := \bC^d$. Note that (\ref{eq_rcom}) take the form $\psi_h^* \psi_k = \delta_{hk} 1$, $\sum_h \psi_h \psi_h^* = 1$, where $\delta_{hk}$ is the Kronecker symbol. 
\end{ex}

\begin{defn} (\cite[\S 1]{DPZ97}) \label{def_hba}
Let $\mA \subset \mB$ be a \sC algebra inclusion. A closed vector space $\mM \subset \mB$ is called {\bf Hilbert $\mA$-bimodule in} $\mB$ if
\begin{enumerate}
\item  $\mM$ is stable w.r.t. left and right multiplication by elements of $\mA$;
\item  $t^* t' \in \mA$, for every $t,t' \in \mM$.
\end{enumerate}
\noindent We say that $\mM$ has {\bf support} $1$ if $\mM \mM^* := {\mathrm{closed \ span}}\left\{ t't^* \ : \  t , t' \in \mM \right\}$ contains an approximate unit for $\mB$.
\end{defn}

Note that if $\mM$ is a Hilbert $\mA$-bimodule in $\mB$, then the map $t ,t' \mapsto t^* t'$ can be regarded as an $\mA$-valued scalar product; moreover, there is a natural identification $\mM \mM^* \simeq K(\mM)$, $t't^* \mapsto \theta_{t,t'}$. The following proposition is a consequence of the universality of the CP-algebra (see \cite[Thm.3.12]{Pim93}).

\begin{prop}
\label{prop_uni_cp}
Covariant morphisms between Hilbert bimodules give rise to \sC algebra morphisms between the associated CP-algebras. In particular, $\com$ is isomorphic to $\mO_{\mM_\alpha}$ for every pullback bimodule $\mM_\alpha$ (Ex.\ref{ex_pullback}). If $\mB$ is a unital \sC algebra, and $\mM$ is a Hilbert $\mA$-bimodule in $\mB$ with support $1$, then there is a canonical morphism $\com \ra \mB$.
\end{prop}

Let $X$ be a locally compact Hausorff space, $\mE \ra X$ a rank $d$ vector bundle, $d \in \bN$. Moreover, let $\wE$ be the Hilbert $C_0(X)$-bimodule of continuous, vanishing at infinity sections of $\mE$, endowed with coinciding left and right $C_0(X)$-module actions. We denote by $\coe$ the CP-algebra associated with $\wE$. For $X$ compact (so that $C(X)$ is unital and $\wE$ is finitely generated), (\ref{eq_rcom}) take the form
\[
\psi_l^* \psi_m = \left \langle \psi_l , \psi_m \right \rangle
\ \ , \ \ 
f \psi_l = \psi_l f 
\ \ , \ \ 
\sum_l \psi_l \psi_l^* = 1 \ \ ,
\]
\noindent $f \in C(X)$. It is proved in \cite[Prop.4.2]{Vas04} that $\coe$ is a locally trivial continuous bundle over $X$, with fibre the Cuntz algebra $\mO_d$. Moreover, it is clear that the circle action (\ref{def_ca}) is by $C_0(X)$-automorphisms: $\alpha : \bT \ra {\bf aut}_X \coe$. In particular, if $\mL \ra X$ is a line bundle, then the fibre of $\mO_\mL$ is the \sC algebra $C(S^1)$ of continuous functions over the circle; CP-algebras associated with line bundles have been classified in \cite[Prop.4.3]{Vas}. In the rest of the present paper, we will consider only vector bundles with rank $> 1$.

The main motivation of the present paper is the classification of the \sC algebras $\coe$ in terms of topological properties of the underlying vector bundles.

\section{Representable $KK$-theory.}
\label{kk_groups}

Let ${\bf C^* alg}$ denote the category of \sC algebras, ${\bf Ab}$ the category of abelian groups. Kasparov constructed a bifunctor $KK_0 :$ ${\bf C^* alg} \times {\bf C^* alg} \ra {\bf Ab}$, assigning to the pair $(\mA , \mB)$ an abelian group $KK_0 (\mA , \mB)$. $KK_0$ depends contravariantly on the first variable, and covariantly on the second one. Let $K_0(\mA)$ denote the $K$-theory of $\mA$, $K^0(\mA)$ the $K$-homology; it turns out that there are isomorphisms $KK (\bC , \mA) \simeq K_0(\mA)$, $KK (\mA , \bC) \simeq K^0(\mA)$.

Let $X$ be a locally compact Hausdorff space, $\mA$, $\mB$ $C_0(X)$-algebras. A $C_0(X)$-{\em Hilbert} $\mA$-$\mB$-{\em bimodule} is a Hilbert $\mA$-$\mB$-bimodule $\mM$ such that $(af) \psi b = a \psi (f b)$ for every $f \in C_0(X)$, $\psi \in \mM$, $a \in \mA$, $b \in \mB$.

Roughly speaking, a $C_0(X)$-Hilbert $\mA$-$\mB$-bimodule can be regarded as the space of sections of a bundle, having as fibres Hilbert $\mA_x$-$\mB_x$-bimodules, $x \in X$.

\begin{ex}
Let $\mA$ be a \sC algebra, $X$ a locally compact Hausdorff space, $\mE \ra X$ a vector $\mA$-bundle in the sense of Mishchenko (\cite{MF80}). Then, the module of continuous, vanishing at infinity sections of $\mE$ has an obvious structure of $C_0(X)$-Hilbert $C_0(X)$-$(X \mA)$-bimodule. 
\end{ex}

\begin{rem} 
\label{rem_tp} 
In the sequel, we will make use of the following two notions of tensor product of $C_0(X)$-Hilbert bimodules.
\begin{enumerate}
\item Let $\mM$ be a $C_0(X)$-Hilbert $\mA$-$\mB$-bimodule, $\mN$ a $C_0(X)$-Hilbert $\mB$-$\mC$-bimodule. We consider the algebraic tensor product $\mM \odot_\mB \mN$ with coefficients in $\mB$, and denote by $\mM \otimes_\mB \mN$ the completition w.r.t. the $\mC$-valued scalar product $\left \langle \psi \otimes \varphi , \psi' \otimes \varphi' \right \rangle :=$ $\left \langle \varphi ,  \left \langle \psi , \psi' \right \rangle_\mM \varphi' \right \rangle_\mN$, $\psi , \psi ' \in \mM$, $\varphi$, $\varphi' \in \mN$; here $\left \langle \cdot , \cdot \right \rangle_\mN$ (resp.  $\left \langle \cdot , \cdot \right \rangle_\mM$) denotes the scalar product on $\mN$ (resp, $\mM)$; note that $\left \langle \psi , \psi' \right \rangle_\mM \in$ $\mB$, so that it makes sense to consider $\left \langle \psi , \psi' \right \rangle_\mM \varphi'$. $\mM \otimes_\mB \mN$ is a $C_0(X)$-Hilbert $\mA$-$\mC$-bimodule in a natural way, and is called the {\bf internal tensor product} of $\mM$ and $\mN$.
\item Let $\mM'$ be a $C_0(X)$-Hilbert $\mA'$-$\mB'$-bimodule. The algebraic tensor product  $\mM \odot_{C_0(X)} \mM'$ with coefficients in $C_0(X)$ is endowed with a natural left $(\mA \otimes_X \mA')$-module action, and with natural $(\mB \otimes_X \mB' )$-valued scalar product and right action. The corresponding completition $\mM \otimes_X \mM'$ is a $C_0(X)$-Hilbert $(\mA \otimes_X \mA')$-$(\mB \otimes_X \mB' )$-bimodule, and is called the {\bf external tensor product} of $\mM$ and $\mM'$.
\end{enumerate} 
\end{rem}

Let $X$ be a $\sigma$-compact metrisable space. Motivated by the Novikov conjecture, Kasparov generalized the construction of $KK_0 ( - , - )$ to the category of $C_0(X)$-algebras (\cite[2.19]{Kas88}); the corresponding bifunctor is called {\em representable} $KK$-theory. We will denote it by the notation $KK ( X ; -,-)$ (note that in \cite{Kas88} the notation $\mR KK (X ; - , -)$ is used). The rest of the present section is devoted to a brief exposition of the construction of $KK (X ; - ,-)$.

\begin{defn}
Let $\mA$, $\mB$ be separable $C_0(X)$-algebras. A {\bf Kasparov $\mA$-$\mB$-module} is a pair $(\mM , F)$, where $\mM$ is a countably generated $C_0(X)$-Hilbert $\mA$-$\mB$-bimodule, and $F = F^* \in L (\mM)$ is an operator such that $[F,a] , a(F^2 - 1) \in K (\mM)$ for every $a \in \mA$. We denote by $\bE (X ; \mA ,\mB)$ the set of Kasparov $\mA$-$\mB$-modules.
\end{defn}

It is customary to consider a $\bZ_2$-grading on Kasparov modules (\cite[\S 14]{Bla}). Since we do not need such a structure, we will assume that every \sC algebra (Hilbert bimodule) is endowed with the trivial $\bZ_2$-grading.

\begin{ex} 
\label{ex_picard}
Let $\mM$ be a countably generated $C_0(X)$-Hilbert $\mA$-$\mB$-bimodule. Let us suppose that $a \in K(\mM)$ for every $a \in \mA$; then $(\mM , 0) \in \bE (X ; \mA , \mB)$, where $0$ is the zero operator. In particular, if every element of $K(\mM)$ is the image of some element of $\mA$ w.r.t. the left $\mA$-module action (so that, there is an isomorphism $\mA \simeq K(\mM)$), then $\mM$ is called {\bf imprimitivity} $\mA$-$\mB$-{\bf bimodule} (see \cite{BGR77}). 
\end{ex}

\begin{ex} 
\label{ex_mor}
Let $\phi : \mA \ra \mB$ be a nondegenerate $C_0(X)$-algebra morphism. We endow $\mB$ with the $C_0(X)$-Hilbert $\mA$-$\mB$-bimodule structure
\[
a , \psi \mapsto \phi (a) \psi      \ \ , \ \
\psi , b \mapsto \psi b     \ \ , \ \
\left \langle \psi , \psi' \right \rangle  := \psi^* \psi' \ ,
\]
\noindent $a \in \mA$, $b , \psi , \psi' \in \mB$, and denote by $\mB_\phi$ the associated $C_0(X)$-Hilbert $\mA$-$\mB$-bimodule. Now, it is clear that $K(\mB_\phi) \simeq \mB$; thus, $\phi (a) \in K(\mB_\phi)$ for every $a \in \mA$. If $\mB$ is $\sigma$-unital (\cite[12.3]{Bla}), then $(\mB_\phi , 0) \in \bE(\mA , \mB)$ (in fact, $\mB_\phi$ is countably generated by an approximate unit $\left\{ u_n \right\}_{n \in \bN} \subset \mB$). If $\mA = \mB$ and $\phi \in {\bf aut}_X \mA$, then the CP-algebra $\mO_{\mA_\phi}$ is isomorphic to the crossed product $\mA \rtimes_\phi \bZ$ (\cite[\S 2]{Pim93}). If $\iota : \mA \ra \mA$ is the identity automorphism, we define $[1]_\mA := (\mA_\iota , 0) \in \bE ( X ; \mA , \mA)$.  
\end{ex}

There are natural notions of {\em homotopy } and {\em direct sum } over $\bE ( X ; \mA ,\mB)$. The {\em representable} $KK${\em -theory group} $KK ( X ; \mA , \mB )$ is constructed in the same way as the usual $KK_0$-group, by endowing the set of homotopy classes of Kasparov $\mA$-$\mB$-modules with the operation of direct sum. The bifunctor $KK(X; - , -)$ is {\em stable}, i.e. $KK ( X ; \mA , \mB )$ is invariant by tensoring $\mA$ or $\mB$ by the \sC algebra $\mK$ of compact operators over a separable Hilbert space.

With an abuse of the notation, we will identify the elements of $\bE ( X ; \mA , \mB )$ with the corresponding classes in $KK ( X ; \mA , \mB )$. We use the notation $+$ to denote the operation of direct sum in $KK ( X ; \mA , \mB )$. When $X = \bullet$ reduces to a single point, then $KK ( \bullet , \mA , \mB )$ is the usual $KK$-group $KK_0 (\mA , \mB)$.

Let $\mA$, $\mB$, $\mA'$, $\mB'$, $\mC$ be separable $C_0(X)$-algebras. We recall that the {\em Kasparov product} (\cite[\S 2.21]{Kas88}) induces bilinear maps
\[
\left\{
\begin{array}{ll}
\times_\mB : KK ( X ; \mA , \mB) \otimes KK ( X ; \mB , \mC) \ra KK ( X ; \mA , \mC) \ ,
\\
\times : KK ( X ; \mA , \mB) \otimes KK ( X ; \mA' , \mB') \ra KK ( X ; \mA \otimes_X \mA'  , \mB \otimes_X \mB') \ .
\end{array}
\right.
\]

\begin{rem}
\label{rem_kas}
In some particular cases, the Kasparov product takes a simple form; in fact, with the notation of Ex.\ref{ex_picard}, Ex.\ref{ex_mor}, we find that

\begin{enumerate}
\item  $(\mM , 0) \times_\mB (\mN , 0) =$ $(\mM \otimes_\mB \mN , 0)$, where $(\mM , 0) \in$ $KK (X ; \mA , \mB)$, $(\mN , 0) \in$ $KK (X ; \mB , \mC)$, and $\otimes_\mB$ denotes the internal tensor product (Rem.\ref{rem_tp});
\item  $(\mM , 0) \times (\mM' , 0) = (\mM \otimes_X \mM' , 0)$, where $(\mM , 0) \in$ $KK (X ; \mA , \mB)$, $(\mM' , 0) \in$ $KK (X ; \mA' , \mB')$, and $\otimes_X$ denotes the external tensor product (Rem.\ref{rem_tp});
\item  Let $\phi : \mA \ra \mB$, $\eta : \mB \ra \mC$ be $C_0(X)$-algebra morphisms. Then $(\mB_\phi , 0) \times_\mB ( \mC_\eta , 0 ) = ( \mC_{\eta \circ \phi} , 0)$.
\end{enumerate} 
\end{rem}

\noindent Thus, $( KK ( X ; \mA , \mA) \ , \ + \ , \ \times_\mA )$ is a ring with identity the class $[1]_\mA$ defined in Ex.\ref{ex_mor}. Let us now consider the ring
\[  RK^0 (X) := KK ( X ; C_0(X) , C_0(X) ) \ \ ; \]
\noindent it is proven in \cite[2.19]{Kas88} that $( RK^0(X) , + )$ is isomorphic to the representable $K$-theory group introduced by Segal in \cite{Seg70}. If $X$ is compact, it is verified that $( RK^0(X) , + )$ coincides with the topological $K$-theory $K^0(X)$. In order for a more concise notation, we denote by $[1]_X \in RK^0(X)$ the class $( [C_0(X)]_\iota ,0 )$ defined in Ex.\ref{ex_mor}.

\begin{lem}
\label{lem_vb_kk}
Let $X$ be a $\sigma$-compact Hausdorff space, $d \in \bN$, $\mE \ra X$ a rank $d$ vector bundle. Then, $C_0(X)$ acts on the left on $\wE$ by elements of $K (\wE)$, and the pair $(\wE , 0)$ is a Kasparov module with class $[\mE] := (\wE , 0) \in RK^0(X)$.
\end{lem}

\begin{proof}
Let $1$ be the identity over $\wE$, $\theta_{\psi,\psi'} \in K (\wE)$, $\psi , \psi' \in \wE$, the operator defined in (\ref{def_theta}). $X$ being $\sigma$-compact, there is a sequence $\left\{ K_n \right\}_n$ of compact subsets covering $X$. Let $\left\{ \lambda_n \right\}$ be a partition of unity with $\ovl{\mathrm{supp} \lambda_n} = K_n$, $n \in \bN$. By the Serre-Swan theorem, the bimodule of continuous sections of the restriction $\left. \mE \right|_{K_n}$ is finitely generated by a set $\left\{ \right. \varphi_{n,k} \left. \right\}_k$; we define $ \psi_{n,k} := \lambda_n \varphi_{n,k} \in \wE$. Let now $u_n := \sum_k \theta_{ { \psi_{n,k}} , { \psi_{n,k}} } \in K (\wE)$. Note that $u_n = \lambda_n^2 \sum_k \theta_{ {\varphi_{n,k}} , {\varphi_{n,k}} } = \lambda_n^2$. Thus, the sequence $U_m := \sum_n^m u_n = \sum_n^m \lambda^2_n$ converges to $1$ in the strict topology. We conclude that $\wE$ is countably generated as a right Hilbert $C_0(X)$-module by the set $\left\{ \psi_{n,k} \right\}$. Let now $f \in C_0 (X)$. We regard at $f$ as an element of $L(\wE)$. Now, $\left\| f - f \sum_n^m u_n \right\| = \left\| f - \sum_n^m \lambda_n^2 f \right\| \stackrel{m}{\ra} 0$; thus $f$ is norm limit of elements of $K (\wE)$, and $C_0(X)$ acts on the left over $\wE$ by elements of $K (\wE)$. We conclude that the pair $(\wE , 0)$ defines a class in $RK^0 (X)$.
\end{proof}

\bigskip

Let $\mA$ be a $C_0(X)$-algebra, $\mE \ra X$ a vector bundle. We define $\wE \otimes_X \mA$ as the external tensor product $\wE \otimes_X \mA_\iota$ (where $\mA_\iota$ is defined in Ex.\ref{ex_mor}). The Kasparov product induces a natural structure of $RK^0(X)$-bimodule on $KK (X ; \mA , \mB)$ (\cite[2.19]{Kas88}). In particular, there is a morphism of unital rings
\begin{equation}
\label{eq_stru_mor}
i_\mA : RK^0(X) \ra KK ( X ; \mA , \mA ) 
\ \ , \ \ 
i_\mA (\mM , F) := (\mM , F) \times [1]_\mA
\ .
\end{equation}
\noindent If $\mE \ra X$ is a vector bundle, then it turns out that $i_\mA [\mE] = ( \wE \otimes_X \mA , 0 )$.

\bigskip

Let $\mA$ be a $C_0(X)$-algebra, $\pic$ the set of isomorphism classes of imprimitivity $C_0(X)$-Hilbert $\mA$-bimodules (Ex.\ref{ex_picard}). We endow $\pic$ with the operation of internal tensor product $\otimes_\mA$; note that the bimodule $\mA_\iota$ defined in Ex.\ref{ex_mor} is a unit for $\pic$. By applying the argument of \cite[\S 3]{BGR77}, it is verified that if $\mM$ is an imprimitivity $C_0(X)$-Hilbert $\mA$-bimodule, and $\ovl \mM$ is the conjugate bimodule, then $\ovl \mM$ is an imprimitivity bimodule, and the map
\[
\mM \otimes_\mA \ovl \mM \ra K(\mM) \simeq \mA_\iota \ \ , \ \ 
\psi' \otimes \ovl \psi \mapsto \theta_{\psi',\psi}
\]
\noindent defines an isomorphism of $C_0(X)$-Hilbert $\mA$-bimodules ($\theta_{\psi',\psi} \in K(\mM)$ is defined by (\ref{def_theta})). Thus, $\pic$ is a group, called the {\em Picard group} of $\mA$.

Let ${\bf {out}}_X \mA$ denote the group of $C_0(X)$-automorphisms of $\mA$ modulo inner automorphisms induced by unitaries in $M(\mA)$. If $\mA$ is $\sigma$-unital, then by \cite[Cor.3.5]{BGR77} we obtain a group anti-isomorphism
\begin{equation}
\label{iso_pic}
\theta : {\mathrm{Pic}} ( X ; \mA \otimes \mK) 
\stackrel{\simeq}{\longrightarrow}
{\bf {out}}_X (\mA \otimes \mK) \ \ , \ \  \mM \mapsto \theta_\mM \ .
\end{equation}
\noindent The previous isomorphism has to be intended in the sense that every imprimitivity $(\mA \otimes \mK)$-bimodule is isomorphic to a bimodule of the type described in Ex.\ref{ex_mor}.

\begin{ex} 
\label{ex_picx}
Let $X$ be a paracompact Hausdorff space. Then ${\mathrm{Pic}}(X ; C_0(X))$ is isomorphic to the Cech cohomology group $H^2(X , \bZ)$. In fact, imprimitivity $C_0(X)$-Hilbert $C_0(X)$-bimodules correspond to modules of continuous sections of line bundles over $X$; on the other hand, it is well-known that the set of line bundles, endowed with the operation of tensor product, is isomorphic as a group to $H^2(X , \bZ)$ (see \cite[\S 3]{BGR77}). 
\end{ex}

We denote by $KK ( X ; \mA , \mA)^{-1}$ the multiplicative group of invertible elements of $KK ( X ; \mA , \mA)$. The argument of Rem.\ref{rem_kas} implies that there is a group morphism
\begin{equation}
\label{eq_mor_pic}
\pi_\mA : \pic  \ra KK ( X ; \mA , \mA)^{-1}  \ \ , \ \
\mM \mapsto (\mM , 0)
\ \ .
\end{equation}

\begin{ex}
We refer to Ex.\ref{ex_picx}. Let $X$ be a locally compact, paracompact Hausdorff space. Then, we have a group morphism
\[
\pi_X : H^2( X , \bZ ) \ra RK^0 (X)^{-1} \ \ ,
\]
\noindent assigning to the isomorphism class of a line bundle the corresponding class in $K$-theory. 
\end{ex}

\section{Graded $\mA$-bundles.}
\label{grading}

Aim of the present section is to assign $KK$-theoretical invariants to $\mA$-bundles carrying a suitable circle action.

\subsection{Circle actions.}

Let $\mA$ be a \sC algebra carrying an automorphic action $\alpha : \bT \ra {\bf aut} \mA$. The \sC dynamical system $(\mA , \bT)$ is said {\em full } if $\mA$ is generated as a \sC algebra by the disjoint union of the {\em spectral subspaces } 
\[
\mA^k := 
\left\{  a \in \mA : \alpha_z (a) = z^k a \ , \ z \in \bT \right\} \ \ , \ \ 
k \in \bZ \ .
\]
\noindent Note that $\mA^h \cdot \mA^k \subseteq \mA^{h+k}$, $(\mA^k)^* = \mA^{-k}$, $h,k \in \bZ$. In particular, every $\mA^k$ is a {\em Hilbert} $\mA^0$-{\em bimodule in} $\mA$ (Def.\ref{def_hba}). Note that $\mA^k$ is full as a right Hilbert $\mA^0$-module if and only if $\mA^{-k} \cdot \mA^k = \mA^0$. Moreover, there is a natural map $\mA^k \cdot \mA^{-k} \ra K(\mA^k)$, $t' t^* \mapsto \theta_{t',t}$.

A \sC dynamical system $(\mA , \bT )$ is said {\em semi-saturated} if $\mA$ is generated as a \sC algebra by $\mA^0$, $\mA^1$ (see \cite{Exe94,AEE95}). It is clear that if $\mA$ is semi-saturated, then $\mA$ is full. From the above considerations, we obtain the following lemma.

\begin{lem}
Suppose $\mA^0 = \mA^1 \cdot \mA^{-1} = \mA^{-1} \cdot \mA^1$. Then, $\mA^1$ is an imprimitivity $\mA^0$-bimodule; if $\mA^1$ is countably generated, the class $\delta_1 (\mA) := (\mA^1 , 0) \in KK_0 ( \mA^0 , \mA^0)$ is defined.
\end{lem}

\begin{ex} 
\label{ex_cp}
Let $\mM$ be a full Hilbert $\mA$-bimodule, $\com$ the associated CP-algebra. Then, $\com$ is semi-saturated w.r.t. the circle action (\ref{def_ca}), so that every $\com^k$, $k \in \bZ$, is an imprimitivity bimodule over the zero-grade algebra $\com^0$. 
\end{ex}

The previous example is universal, as we will show in the next lemma. 

Let us introduce the following terminology: if $(\mA , \bT , \alpha)$, $(\mB , \bT , \beta)$ are \sC dynamical systems, a {\em graded morphism} is a \sC algebra morphism $\phi : \mA \ra \mB$ such that $\phi (\mA^k) \subseteq \mB^k$, $k \in \bZ$. Graded morphisms will be denoted by the notation $\phi : (\mA , \bZ) \ra (\mB , \bZ)$.

\begin{lem}(\cite[Thm.3.1]{AEE95}) \label{lem_isop}
Let $(\mA , \bT )$ be semi-saturated, and $\mA^1$ full as a Hilbert $\mA^0$-bimodule. Then, there is a graded isomorphism $(\mA , \bZ) \simeq ( \mO_{\mA^1} , \bZ )$, where $\mO_{\mA^1}$ is the CP-algebra associated with $\mA^1$.
\end{lem}

\begin{proof}
It is a direct consequence of Prop.\ref{prop_uni_cp}: in fact, $\mA^1$ is a Hilbert $\mA^0$-bimodule in $\mA$ with support $1$, and generates $\mA$ as a \sC algebra. 
\end{proof}

\subsection{Graded Bundles.}

As usual, we denote by $X$ a locally compact Hausdorff space.

\begin{defn}
Let $(\mA , G , \alpha )$ be a \sC dynamical system. An $\mA$-bundle $(\mF , (\pi_x : \mF \ra \mA)_{x \in X} )$ has a {\bf global $G$-action} if there exists an action $\alpha^X : G \ra {\bf aut}_X \mF$, such that $\pi_x \circ \alpha^X = \alpha \circ \pi_x$ for every $x \in X$.
\end{defn}

Let $\mF$ be a $\mA$-bundle carrying a global $\bT$-action. Then, every spectral subspace $\mF^k$, $k \in \bZ$ has a additional structure of $C_0(X)$-Hilbert $\mF^0$-bimodule, in fact $f t = t f$, $t \in \mF^k$, $f \in C_0(X)$. We say in such a case that $\mF$ is a {\em graded} $\mA$-bundle. 

\begin{prop}
\label{thm_f1}
Let $(\mA , \bT , \alpha)$ be a semi-saturated \sC dynamical system, with $\mA^1$ full as a Hilbert $\mA^0$-bimodule. Moreover, let $X$ be paracompact. Then, for every graded $\mA$-bundle $\mF$ over $X$ there is an isomorphism $(\mF , \bZ) \simeq (\mO_{\mF^1} , \bZ)$. If $\mB$ is a graded $\mA$-bundle, there is an isomorphism $(\mF , \bZ) \simeq (\mB , \bZ)$ if and only if the Hilbert bimodules $\mF^1$, $\mB^1$ are covariantly isomorphic.
\end{prop}

\begin{proof}
We prove that $( \mF , \bT , \alpha^X )$ is semi-saturated, and that $\mF^1$ is full as a Hilbert $\mF^0$-bimodule. Let $\mU := \left\{ U \subseteq X \right\}$ be an open (locally finite) trivializing cover (i.e., every restriction $\mF_U := C_0(U) \mF$ is isomorphic to $C_0(U) \otimes \mA$, $U \in \mU$). Since $\alpha^X (t) \in \mF_U$, $t \in \mF_U$, for every $U \in \mU$ we obtain a global action $\alpha_U : \bT \ra {\bf aut}_U \mF_U$. Since $\mF_U$ is a trivial bundle, it is clear that $( \mF_U  , \bT , \alpha^U )$ is semi-saturated, and that $\mF^1_U$ is full as a Hilbert $\mF^0_U$-bimodule. We now consider a partition of unit $\left\{ \lambda_U \in C_0(X) \right\}$ subordinate to $\mU$; since every $t \in \mF$ admits a decomposition $t = \sum_U \lambda_U t$, with $\lambda_U t \in \mF_U$, we conclude that $( \mF , \bT , \alpha^X )$ is semi-saturated, and that $\mF^1$ is full. By applying Lemma \ref{lem_isop}, we obtain the isomorphism $(\mF , \bZ) \simeq (\mO_{\mF^1} , \bZ)$. The second assertion is an immediate consequence of the $\bZ$-grading defined on $\mF$, $\mB$.
\end{proof}

\begin{cor}
\label{morita}
With the above notation, every $\mF^r$ is an imprimitivity $\mF^0$-bimodule, $r \in \bN$.
\end{cor}

\begin{proof}
It suffices to consider the identifications $K (\mF^r) \simeq \mF^r \cdot (\mF^r)^* \simeq \mF^0$.
\end{proof}

%
%
%
%
%

From Prop.\ref{thm_f1}, we have an interpretation of the set of isomorphism classes of graded $\mA$-bundles in terms of covariant isomorphism classes of $C_0(X)$-Hilbert bimodules. Thus, a description in terms of $KK(X ; - , -)$-groups becomes natural. As a first step, we consider the zero grade algebra; for every \sC algebra $\mA$, we denote by $H^1(X , {\bf aut}\mA)$ the set of isomorphism classes of $\mA$-bundles over $X$ (see for example \cite[Thm.2.1]{Vas} for a justification of such a notation). $H^1(X , {\bf aut}\mA)$ has a distinguished element, called $0$, corresponding to the trivial $\mA$-bundle.

\bigskip

Let $\mF$ be a graded $\mA$-bundle over $X$. We define
\begin{equation}
\label{def_delta_0}
\delta_0 (\mF) := [ \mF^0 \otimes \mK ] \in H^1(X,{\bf aut} ( \mA^0 \otimes \mK) ) \ \ .
\end{equation}
\noindent Thus, the equality $\delta_0 (\mF) = \delta_0 (\mB)$ is intended in the sense that $\mF^0$, $\mB^0$ are stably isomorphic as $\mA^0$-bundles.

\begin{rem}
\label{nistor}
Let $X$ be a pointed, compact, connected $CW$-complex such that the pair $(X , x_0)$, $x_0 \in X$, is a homotopy-cogroup. We denote by $SX$ the (reduced) suspension. In order for more compact notations, we define $X^\bullet := X - \left\{ x_0 \right\}$. It follows from a result by Nistor (\cite[\S 5]{Nis92}) that 
\[
H^1 (SX , {\bf aut}\mO_d^0) \simeq   [X , {\bf aut} \mO_d^0 ]_{x_0}  
                            \simeq   KK_1 ( C \mO_d^0 \ , \ X^\bullet \mO_d^0 ) \ ,
\]
\noindent where $C \mO_d^0 := \left\{ ( z , a) \in \bC \oplus C_0 ( [0,1) , \mO_d^0 ) \ : \ a (0) = z 1 \right\}$ is the mapping cone. We also find
\[
H^1( SX ,  {\bf aut} (\mO_d^0 \otimes \mK) ) =
[X , {\bf aut} (\mO_d^0 \otimes \mK) ]_{x_0}  
\simeq   KK_0 ( \mO_d^0 \ , \ X^\bullet \mO_d^0 ) \ .
\]
\noindent In particular, when $X$ is the $(n-1)$-sphere, we obtain $H^1 ( S^n , {\bf aut} \mO_d^0  ) = \left\{ 0 \right\}$, as proved also in \cite[Thm.1.15]{Tho87}.
\end{rem}

\

Let $X$ be a $\sigma$-compact metrisable space, $\mA$ separable and $\sigma$-unital; then, every graded $\mA$-bundle $\mF$ over $X$ is separable and $\sigma$-unital, and $\mF^1$ is countably generated as a Hilbert $\mF^0$-bimodule (in fact, $\mF$ is countably generated over compact subsets). Moreover, Cor.\ref{morita} implies that $\mF^0$ acts on the left on $\mF^1$ by elements of $K (\mF^1)$. Thus, we define
\begin{equation}
\label{def_delta_1}
\delta_1 (\mF) := ( \mF^1 ,0 )  \in KK ( X ; \mF^0 , \mF^0 ) \ \ .
\end{equation}
\noindent With an abuse of notation, in the sequel we will denote by $\delta_1 (\mF)$ also the class of $(\mF^1 , 0)$ in $KK_0 ( \mF^0 , \mF^0 )$ obtained by {\em forgetting} the $C_0(X)$-structure. Note that since $\mF^1$ is an imprimitivity bimodule, we find that $\delta_1 (\mF)$ is invertible; thus, the Kasparov product by $\delta_1 (\mF)$ defines an automorphism on $KK_0 ( \mB , \mF^0 )$ for every \sC algebra $\mB$. In particular, $\delta_1 (\mF) \in {\bf aut} K_0 (\mF^0)$.

From \cite[Thm.4.9]{Pim93} and Prop.\ref{thm_f1}, we get an exact sequence for the $KK$-theory of $\mF$. It is clear that in the case in which $\mF$ is the CP-algebra of a vector bundle, we may directly apply \cite[Thm.4.9]{Pim93} by replacing $\mF^0$ with $C_0(X)$.

\begin{cor}
\label{cor_ex_seq}
For every separable \sC algebra $\mB$, and graded $\mA$-bundle $\mF$, the following exact sequence holds:
\begin{equation}
\label{ex_seq}
\xymatrix{
           KK_0 (\mB , \mF^0)
		    \ar[r]^-{{1 - \delta_1(\mF)}}
		 &  KK_0 (\mB , \mF^0)
		    \ar[r]^-{i_0}
		 &  KK_0 (\mB , \mF)
		    \ar[d]^-{\delta_0}
		 \\ KK_1 (\mB , \mF)
		    \ar[u]^-{\delta_1}
		 &  KK_1 (\mB , \mF^0)
		    \ar[l]_-{i_1 }
		 &  KK_1 (\mB , \mF^0)                             
		    \ar[l]_{{1- \delta_1(\mF)}}
}
\end{equation}
\noindent where $i_*$ are the morphisms induced by the inclusion $\mF^0 \hra \mF$, and $\delta_*$ are the connecting maps induced by the $KK$-equivalence between $\mF^0$, $\mT_{\mF^1}$.
\end{cor}

We introduce a notation. Let $\mF$, $\mB$ be graded $\mA$-bundles with $\delta_0(\mF) = \delta_0(\mB)$; then, there is a $C_0(X)$-algebra isomorphism $\alpha : \mF^0 \otimes \mK \ra \mB^0 \otimes \mK$, and a ring isomorphism $\alpha_* : KK (X ; \mF^0 , \mF^0) \ra KK (X ; \mB^0 , \mB^0)$ is defined. We write

\begin{equation}
\label{def_iso_delta}
\delta (\mF) = \delta (\mB)  \ \
\Leftrightarrow \ \ 
\delta_0 (\mF) = \delta_0 (\mB) \ \
{\mathrm {and}} \ \
\alpha_* \delta_1(\mF) = \delta_1 (\mB) \ ;
\end{equation}

\noindent note that we used the stability of $KK(X ; - , - )$, so that we identified $( \mF^1 , 0 ) \in KK (X ; \mF^0 , \mF^0)$ with $( \mF^1 \otimes \mK , 0) \in KK (X ; \mF^0 \otimes \mK , \mF^0 \otimes \mK )$. The tensor product of $\mF^1$ by $\mK$ is understood as the external tensor product of Hilbert bimodules. Note that $\alpha_* \delta_1(\mF) =  ( (\mF^1 \otimes \mK)_\alpha , 0 )$, where $(\mF^1 \otimes \mK)_\alpha$ is the pullback bimodule defined as in Ex.\ref{ex_pullback}. Let us now consider the natural $\bZ$-gradings on $\mF \otimes \mK$, $\mB \otimes \mK$ induced by $\mF$, $\mB$; if there is an isomorphism $\alpha : (\mF \otimes \mK , \bZ) \ra (\mB \otimes \mK , \bZ)$, then $\delta (\mF) = \delta (\mB)$.

\begin{ex}
Let $\mF := X \mA$; then, $\delta_1 (\mF) = \delta_1 (\mA) \times [1]_X$  (we used the notation (\ref{def_xa})). 
\end{ex}

\noindent Let $\mF$ be a graded $\mA$-bundle. In general, it is clear that elements of $KK ( X ; \mF^0 , \mF^0 )$ do not arise from grade-one components of graded $\mA$-bundles. Anyway, they can be recognized by considering any open trivializing cover $\mU := \left\{ U \subseteq X \right\}$ for $\mF$, and by noting that the conditions $\delta_0 ( \mF^0_U ) = 0$, $\delta_1 (\mF^1_U ) = \delta_1 (\mA) \times [1]_U$ hold (see previous example).

\

\begin{rem}
\label{rem_f}
Let $\mF$, $\mB$ be graded $\mA$-bundles over a $\sigma$-compact Hausdorff space $X$, with a fixed $C_0(X)$-isomorphism $\alpha : (\mF^0 \otimes \mK , \bZ ) \ra (\mB^0 \otimes \mK , \bZ )$. The pullback bimodule $(\mF^1 \otimes \mK)_\alpha$ has the same class as $\mB^1 \otimes \mK$ in ${\mathrm{Pic}}(X ; \mB^0 \otimes \mK)$ if and only if there is a $C_0(X)$-isomorphism $(\mF \otimes \mK , \bZ ) \ra (\mB \otimes \mK , \bZ )$ (see Prop.\ref{thm_f1}). Let $\theta(\mF) := \theta_{\mF^1 \otimes \mK} \in {\bf aut}_X (\mF^0 \otimes \mK)$ be defined according to (\ref{iso_pic}), up to inner automorphisms. The CP-algebra associated with $(\mF^0 \otimes \mK)_{\theta(\mF)}$ is graded isomorphic to $(\mF^0 \otimes \mK) \rtimes_{\theta(\mF)} \bZ$ (see Ex.\ref{ex_mor}); by universality, we obtain the isomorphism
\[ 
( \mF \otimes \mK , \bZ )
\simeq 
( ( \mF^0 \otimes \mK ) \rtimes_{\theta(\mF)} \bZ  , \bZ ) \ . 
\]
\noindent Thus, ${\mathrm{Pic}}(X ; \mF^0 \otimes \mK)$ describes the set of graded isomorphism classes of stabilized $\mA$-bundles. The map
\[
{\mathrm{Pic}}(X ; \mF^0)  \ra  KK ( X ; \mF^0 , \mF^0 )  \ \ , \ \ 
\theta (\mF) \mapsto \delta_1 (\mF)
\]
\noindent gives a measure of the accuracy of the class $\delta_1$ in describing the set of graded isomorphism classes of stabilized $\mA$-bundles. 
\end{rem}

\begin{ex} 
Let $\mE \ra X$ be a rank $d$ vector bundle. We denote by $\alpha^X : \bT \ra {\bf aut}_X \coe$ the circle action (\ref{def_ca}). Let now $\pi_x : \coe \ra \mO_d$, $x \in X$, be the fibre epimorphisms of $\coe$ as an $\mO_d$-bundle; we denote by $\alpha_x : \bT \ra {\bf aut} \mO_d$ the circle action (\ref{def_ca}) on the Cuntz algebra. By definition, it turns out that $\pi_x \circ \alpha^X = \alpha_x \circ \pi_x$, thus $\alpha^X$ is a global $\bT$-action. Moreover, $\alpha^X$ is full and semi-saturated, according to Ex.\ref{ex_cp}. Thus, the previous considerations apply with $\mF = \coe$, $\mA = \mO_d$; in particular, $\coe^1$ is an imprimitivity $C_0(X)$-Hilbert $\coe^0$-bimodule. 
\end{ex}

\subsection{The $KK$-class for the CP-algebra of a vector bundle.}

Let $\mE \ra X$ be a rank $d$ vector bundle over a $\sigma$-compact Hausdorff space $X$. We denote by $i : RK^0(X) \ra KK (X ; \coe^0 , \coe^0)$ the structure morphism (\ref{eq_stru_mor}). In order to simplify the notation, we write ${\bf 1}_d := [1]_{\mO_d^0} \in KK_0(\mO_d^0 , \mO_d^0)$. The following result has been proved in \cite[Thm.5.6]{Vas04}.

\begin{thm}
\label{delta1_coe}
With the above notation,
\begin{equation}
\label{eq_main}
\delta_1 (\coe) = i[\mE] = [\mE] \times [1]_{\coe^0}  \ .
\end{equation}
\noindent In particular, $\delta_1 (X \mO_d) = d [1]_X \times {\bf 1}_d$.
\end{thm}

We now discuss the properties of the class $\delta_1 (\coe)$ in the case in which the base space is an even sphere $S^{2n}$; for this purpose, recall that $K^0(S^{2n}) = \bZ^2$, $K^1(S^{2n}) = 0$. We will also make use of the ring structure of $K^0(S^{2n})$, induced by the operation of tensor product: it turns out that there is a ring isomorphism $K^0(S^{2n}) \simeq \bZ^2 [\lambda] / ( \lambda^2 )$, i.e., the elements of $K^0(S^{2n})$ are polynomials of the type $z + \lambda z'$, $z,z' \in \bZ$, with the relation $\lambda^2 = 0$ (\cite[Chp.11]{Hus}). Tensoring by a vector bundle $\mE \ra S^{2n}$ with class $d + \lambda c \in K^0(S^{2n})$ defines an endomorphism $\lambda_\mE \in {\bf end} K^0(S^{2n})$, $\lambda_\mE (z + \lambda z') := dz + \lambda (d z' + cz)$. We also denote by $\iota$: $\iota (z + \lambda z') := z + \lambda z'$ the identity automorphism on $K^0(S^{2n})$. Note that $\iota - \lambda_\mE$ is injective (for $d > 1$).

%
%
%
%

We can now compute the $K$-theory of $\coe$ over even spheres: the exact sequence \cite[Thm.4.9]{Pim93}, and the above considerations, imply
\[
\xymatrix{
           \bZ^2
		    \ar[r]^-{{\iota - {\lambda_\mE}}}
		 &  \bZ^2
		    \ar[r]^-{i_0}
		 &  K_0(\coe)
		    \ar[d]^-{0}
		 \\ K_1 (\coe)
		    \ar[u]^-{0}
		 &  0
		    \ar[l]_-{i_1 }
		 &  0 
		    \ar[l]_{0}
}
\]
\noindent so that, we have the $K$-groups
\[
K_0 (\coe) =  \frac{\bZ^2}{(\iota - \lambda_\mE) \bZ^2 }     \ \ , \ \ 
K_1(\coe) = 0 \ \ .
\]
\noindent In the case in which $d-1$ and $c$ are relatively prime, some elementary computations show that $K_0(\coe) = \bZ_{(d-1)^2}$. This shows that $\coe$ is non-trivial as an $\mO_d$-bundle: in fact, the Kunneth theorem (\cite[\S 23]{Bla}) implies $K_0( S^{2n} \mO_d ) = \bZ_{d-1} \oplus \bZ_{d-1}$.

We now pass to describe the class $\delta_1(\coe)$. For every rank $d$ vector bundle $\mE \ra S^{2n}$ we find $\delta_0 (\coe) = 0$, in fact $\coe^0 \simeq S^{2n} \mO_d^0$ (see Rem.\ref{nistor}); we denote by $$\alpha_\mE : KK ( S^{2n}  ;  \coe^0 , \coe^0 ) \ra KK ( S^{2n} \ ; \ S^{2n} \mO_d^0 , S^{2n} \mO_d^0 )$$ the induced ring isomorphism. The Kunneth theorem implies 
\[
\left\{
\begin{array}{ll}
K_0 (S^{2n} \mO_d^0) = 
K^0(S^{2n}) \otimes K_0 (\mO_d^0)  = 
\bZ^2 \otimes \bZ \left[ \frac{1}{d}  \right]   \\
KK_0 ( S^{2n} \mO_d^0 , S^{2n} \mO_d^0 ) = 
{\bf end} \left( \bZ^2 \otimes \bZ \left[ \frac{1}{d}  \right] \right)
\end{array}
\right.
\]
\noindent where $\bZ \left[ \frac{1}{d}  \right]$ is the group of $d$-adic integers (\cite[10.11.8]{Bla}). We denote by $$(z + \lambda z') \otimes q \ \ , \ \ z,z' \in \bZ \ , \ q \in \bZ \left[ \frac{1}{d}  \right] \ ,$$ the generic "elementary tensor" in $K_0(S^{2n} \mO_d^0)$. Let now $$\delta'_1 (\coe) := \alpha_\mE \circ \delta_1 (\coe) = \alpha_\mE \circ i [\mE] \in KK ( S^{2n} \ ; \ S^{2n}\mO_d^0 \ , \ S^{2n}\mO_d^0 )  \  ;$$ it is clear that $\delta'_1 (\coe) = [\mE] \times [1]_{S^{2n}\mO_d^0}$ is the $KK$-class associated with the bimodule $\wE \otimes_X (S^{2n} \mO_d^0)$. We denote by $\beta_\mE \in {\bf aut} K_0(S^{2n} \mO_d^0)$ the automorphism induced by $\delta'_1 (\coe)$; it follows from the above considerations that $\beta_\mE = \lambda_\mE \times {\bf 1}_d \in {\bf end} \left( \bZ^2 \otimes \bZ \left[ \frac{1}{d}  \right] \right)$. Thus, we find
\begin{equation}
\label{def_aute}
\beta_\mE ((z + \lambda z') \otimes q) =  
\lambda_\mE (z + \lambda z') \otimes q  =
[dz + \lambda (d z' + cz)] \otimes q
\ \ .
\end{equation}

\begin{thm}
\label{prop_sn}
Let $\mE , \mE' \ra S^{2n}$ be rank $d$ vector bundles. Then, the following are equivalent:
\begin{enumerate}
\item  $\delta'_1 (\coe) = \delta'_1 (\mO_{\mE'}) \ \in \ KK (  S^{2n} \  ; \ S^{2n} \mO_d^0 \  , \  S^{2n} \mO_d^0  )$;
\item  there is a $C( S^{2n})$-isomorphism $(\coe \otimes \mK , \bZ) \ra ( \mO_{\mE'} \otimes \mK , \bZ )$;
\item  $[\mE] = [\mE'] \in K^0(S^{2n})$.
\end{enumerate}
\end{thm}

\begin{proof}
\

1) $\Rightarrow$ 3): we have $\delta'_1 (\coe) = \delta'_1 (\mO_{\mE'})$ if and only if $\beta_\mE = \beta_{\mE'}$; thus, (\ref{def_aute}) implies $\beta_\mE (1 \otimes 1) = (d + \lambda c) \otimes 1 = (d + \lambda c') \otimes 1 = \beta_{\mE'} (1 \otimes 1)$, where $[\mE] := d + \lambda c$, $[\mE'] := d + \lambda c'$. From the equality $(d + \lambda c) \otimes 1 = (d + \lambda c') \otimes 1$, we conclude $[\mE] = [\mE']$ (note that we know {\em a priori} that $\mE$, $\mE'$ have the same rank $d$).
%
%

3) $\Rightarrow$ 2) follows from \cite[Prop.5.10]{Vas04}.

2) $\Rightarrow$ 1) is trivial by definition of $\delta_1$.
\end{proof}

Let $\mE , \mE' \ra S^{2n}$ be rank $d$ vector bundles with classes $[\mE] = d +\lambda c$, $[\mE'] = d + \lambda c'$, such that $d-1$, $c$ and $d-1$,$c'$ are relatively prime; then, $K_0(\coe) = K_0(\mO_{\mE'}) = \bZ_{(d-1)^2}$. If $c \neq c'$, then the previous theorem implies that $\coe$ is not graded stably isomorphic to $\mO_{\mE'}$, and $\delta_1(\coe) \neq \delta_1(\mO_{\mE'})$. This shows that $\delta_1$ is a more detailed invariant w.r.t. the $K$-theory of $\coe$.

\

We conclude with a remark about odd spheres: in this case $K^0(S^{2n+1}) = \bZ$, so that for every rank $d$ vector bundle $\mE \ra S^{2n+1}$ we find $[\mE] = d$. Thus, \cite[Prop.5.10]{Vas04} implies $( \coe \otimes \mK , \bZ ) \simeq ( S^{2n+1} \mO_d \otimes \mK , \bZ   )$.




\end{document}